\documentstyle[11pt]{article}

\begin{document}

\large

\begin{center}

{\bf Some new moment rearrangement invariant spaces; theory and applications.}\\

\vspace{3mm}

{\sc Eugene Ostrovsky, Leonid Sirota}\\

\vspace{3mm}

 {\it Department of Mathematic, HADAS company, \\
56209, Hertzelia Pituach, Galgaley Aplada street, 16, ISRAEL; \\
 e - mail: galo@list.ru }\\
{\it Bar - Ilan University,  59200, Ramat Gan, ISRAEL; \\
 e - mail: sirita@zahav.net.il }\\

\vspace{3mm}

              Abstract.\\
\end{center}

\textwidth = 8 cm
\begin{verse}
\normalsize
\hspace{5mm} In this article we introduce and investigate
\hspace{5mm} some new Banach spaces,
\hspace{5mm} so - called moment spaces, and consider applications to the Fourier series,
\hspace{5mm} singular integral operators, theory of martingales.\\

\vspace{3mm}


\vspace{3mm}


\end{verse}

\vspace{4mm}

\textwidth = 15 cm
{\bf 1 \ Definitions. Simple Properties.}\\

\vspace{3mm}

 Let $ (X, \Sigma,\mu) $ be a measurable space with
non - trivial measure $ \mu: \ \exists A \in \Sigma, \mu(A) \in (0,\mu(X)). $
We will assume that either
$ \mu(X) = 1, $ or $ \mu(X) = \infty $ and that the measure
 $ \mu \  $ is $ \sigma - $ finite and diffuse:
$ \forall A \in \Sigma, 0 < \mu(A) < \infty \ \exists B \subset A, \mu(B) =
\mu(A)/2. $  Define as usually for all the measurable function
$ f: X \to R^1 \ $
$$
|f|_p = \left(\int_X |f(x)|^p \ \mu(dx) \right)^{1/p}, \ p \ge 1;
$$
$ L_p = L(p) = L(p; X,\mu) = \{f, |f|_p < \infty \}. $
 Let $ a = const \ge 1, b = const \in (a,\infty], $ and let $ \psi = \psi(p) $
be some positive continuous on the {\it open} interval $ (a,b) $
 function, such that  there
exists a measurable function $  f: X \to R $ for which
$$
\psi(p) = |f|_p, \ p \in (a,b).
$$
 Note that the function $ p \to p \cdot \log \psi(p), \ p \in (a,b) $
is convex.\par
 The set of all those functions we will denote $ \Psi: \ \Psi = \Psi(a,b) =
\{ \psi(\cdot) \}. $ We can describe all those functions.\\
{\bf Theorem 0.}  {\it Let the measure } $ \mu $ {\it is diffuse.}
 {\it The function } $ \nu(p), \ p \in (a,b) $ {\it belongs to the set }
$ \Psi $ {\it if and only if there exist a two functions } $ \Lambda_1(p), \
\Lambda_2(p), $ {\it such that } $ \nu^p(p) = \Lambda_1(p) + \Lambda_2(p), $
{\it where } $ \Lambda_1(p) $ {\it is absolute monotonic on the interval }
$ (a,b) $ {\it and } $ \Lambda_2(p) $ {\it is relative monotonic on the
interval } $ (a,b): \ \forall k = 0,1,2,\ldots $
$$
\forall p \in (a,b) \ \Rightarrow \Lambda^{(k)}_1(p) \ge 0, \
(-1)^k \Lambda^{(k)}_2(p) \ge 0.
$$
{\bf Proof.}  Let $ \nu(\cdot) \in \Psi, $ then $ \exists f: X \to R, \
\nu^p(p) = $
$$
\int_X |f(x)|^p \ \mu(dx) = \int_X \exp(p \log |f(x)|)) \mu(dx) =
\Lambda_1(p) + \Lambda_2(p),
$$
where
$$
\Lambda_1(p) = \int_{ \{x: |f(x)| \ge 1 \}} \exp(p \log|f(x)|) \ \mu(dx),
 \ \Lambda^{(k)}_1(p) \ge 0; \
$$
$$
\Lambda_2(p) = \int_{ \{x: |f(x)| < 1 \} } \exp(p \log|f(x)|) \ \mu(dx),
\ (-1)^k \Lambda^{(k)}_2(p) \ge 0.
$$

 Inversely, assume that $ \nu^p(p) = \Lambda_1(p) + \Lambda_2(p), \
\Lambda^{(k)}_1(p) \ge 0, (-1)^{(k)}\Lambda^{(k)}(p) \ge 0. $
It follows from Bernstein's theorem that
$$
\Lambda_1(p) = \int_{R} \exp(pt) \mu_1(dt), \ \Lambda_2(p) =
\int_{R} \exp(pt) \mu_2(dt),
$$
where $ \mu_1, \ \mu_2 $ are a Borel measures on the set $ R $
such that $ supp \ \mu_1 \in [0,
\infty), \ supp \ \mu_2 \in (-\infty,0] $ and

$$
\forall p \in (a,b) \ \Rightarrow
\Lambda_1(p) < \infty, \ \Lambda_2(p) < \infty.
$$
Therefore
$$
\nu^p(p) = \int_{-\infty}^{\infty} \exp(pt) (\mu_1(dt) + \mu_2(dt)).
$$
 Since the measure $ \mu $ is diffuse, there exists a (measurable)
function $ \eta: X \to R $ such that
$$
\nu^p(p) = \int_X \exp(p \eta(x)) \ \mu(dx).
$$
Thus, for $ f(x) = \exp(\eta(x)) $ we obtain:
$$
|f|^p_p = \int_X \exp(p \eta(x)) \mu(dx) = \nu^p(p), \ |f|_p = \nu(p).
$$
{\bf Corollary 1.}
 Note that if $ \psi_1(\cdot) \in \Psi, \ \psi_2(\cdot) \in \Psi, $ then
$ \psi_1(\cdot) \cdot \psi_2(\cdot) \in \Psi. $ Indeed, if
$$
\psi_1(p) = |f_1|_p, \ \psi_2(p) = |f_2|_p,
$$
and the functions $ f_1, \ f_2 $ are independent, then
$$
\psi_1(p) \cdot \psi_2(p) = |f_1 \cdot f_2|_p.
$$

 We extend the set $ \Psi $ as follows:
$$
 EX \Psi \stackrel{def}{=} EX \Psi(a,b) = \{ \nu = \nu(p) \} =
$$
$$
\{ \nu: \exists \psi(\cdot) \in \Psi: 0 < \inf_{p \in (a,b)} \psi(p)/\nu(p)
\le \sup_{p \in (a,b)} \psi(p)/\nu(p) < \infty \},
$$
 $$
U\Psi \stackrel{def}{=} U \Psi(a,b) = \{\psi = \psi(p),
\forall p \in (a,b) \ \Rightarrow \psi(p) > 0 \}
$$
and the function $ p \to \psi(p), \ p \in (a,b) $ is continuous. \par
 Hereafter $ a = const \ge 1, \ b \in (a, \infty]. $ \par
{\it Since the case $ \psi(a+0) < \infty, \ \psi(b-0) < \infty $
is trivial for us,
we will assume further that either} $ \psi(a+0) = \infty $ {\it or}
$ \psi(b-0) = \infty, $ {\it or both the cases:} $ \psi(a+0) = \psi(b-0)
= \infty.$ \par
 We define in the case $ b = \infty \ \psi(b-0) = \lim_{p \to \infty}
\psi(p). $ \par
{\bf Definition 1.} Let $ \psi(\cdot) \in U\Psi(a,b). $  The space
$ G(\psi) = G(X,\psi) = G(X,\psi, \mu) = G(X,\psi,\mu, a,b) $
consist on all the measurable functions $ f: X \to R $ with finite norm
$$
||f||G(\psi) \stackrel{def}{=} \sup_{p \in (a,b)} \left[ |f|_p/\psi(p) \right].
$$

 The spaces $ G(\psi), \ \psi \in U\Psi $ are non - trivial: arbitrary bounded
$ \sup_x |f(x)| < \infty $
measurable function $ f: X \to R $ with finite support: $ \mu(supp \ |f|)
< \infty $ belongs to arbitrary space $ G(\psi), \ \forall \psi \in
U \Psi. $ \par
 We denote as usually $ supp \ \psi = \{p: \ \psi(p) < \infty \}. $ \par
 Now we consider a very important for applications examples
 of $ G(\psi) $ spaces.  Let  $ a = const \ge 1,
b = const \in (a, \infty]; \alpha, \beta = const. $ Assume also that at
$ b < \infty \ \min(\alpha,\beta) \ge 0 $ and denote by
$ h $ the (unique) root of equation
$$
(h-a)^{\alpha} = (b-h)^{\beta}, \ a < h < b;
  \ \zeta(p) = \zeta(a,b; \alpha,\beta; p) =
$$
$$
(p-a)^{\alpha}, \ p \in (a,h);
\ \zeta(a,b; \alpha,\beta;p) = (b-p)^{\beta}, \ p \in [h,b);
$$
and in the case $ b = \infty $ assume that
 $ \alpha \ge 0, \beta < 0; $ denote
by $ h $ the (unique) root of equation
 $ (h-a)^{ \alpha} = h^{ \beta}, h > a; $ define in this case
$$
\zeta(p) = \zeta(a,b;\alpha,\beta;p) = (p-a)^{\alpha}, \  p
\in (a,h); \ p \ge h \ \Rightarrow \zeta(p) = p^{\beta}.
$$
Note that at
 $ b = \infty \ \Rightarrow \zeta(p) \asymp (p-a)^{\alpha} \
p^{-\alpha + \beta} \asymp \min \{(p-a)^{\alpha}, p^{\beta} \},
\ p \in (a,\infty) $ and that at
 $ b < \infty \ \Rightarrow \zeta(p) \asymp
(p-a)^{\alpha} (b-p)^{\beta} \asymp \min \{(p-a)^{\alpha}, (b-p)^{\beta} \},
 p \in (a,b). $  Here and further
 $ p \in (a,b) \ \Rightarrow \psi(p) \asymp \nu(p) $  denotes that
$$
0 < \inf_{p \in (a,b)} \psi(p)/\nu(p) \le \sup_{p \in (a,b)} \psi(p)/\nu(p)
< \infty.
$$

 We will denote also by the symbols $ C_j, j \ge 1 $ some "constructive" finite
non - essentially positive constants. By definition,
 $ I(A) = I(A,x) = I(x \in A) = 1, x \in A; \ I(A) = 0,  x \notin A. $ \\
{\bf Definition 2.} The space $ G = G_X = G_X(a,b;\alpha,\beta)=
G(a,b; \alpha,\beta) $  consists on all measurable functions
 $ f: X \to R^1 $ with finite norm
$$
||f||G(a,b; \alpha,\beta) = \sup_{p \in (a,b)} \left[ |f|_p \cdot
\zeta(a,b; \alpha,\beta;p) \right].
$$
{\bf Corollary 2. } As long as the cases $ \alpha \le 0; \ b < \infty,
\beta \le 0 $ and $ \ b = \infty, \beta \ge 0 $ are trivial,
we will assume further that
either $ 1 \le a < b < \infty, \min(\alpha,\beta) > 0, $ or
$ 1 \le a, b = \infty, \alpha \ge 0, \beta < 0. $ \\
{\bf Lemma 1. }  {\it Let} $ \psi \in U \Psi, \ \psi(a+0) = \psi(b-0) =
\infty, b < \infty. $ {\it There exist a two functions $ \nu_1, \nu_2 \in
U\Psi, \nu_1(a+0) \in (0,\infty),  \nu_1(p) \sim \psi(p), p \to b-0;
\nu_2(b-0) \in (0,\infty), \nu_2(p) \sim \psi(p), p \to a+0 $ {\it such that }
 {\it the space } $ G(\psi) $ may be represented as a direct sum }
$$
G(\psi) = G(\nu_1) + G(\nu_2).
$$
{\bf Proof.}  Indeed, if
 $ f = f_1 + f_2, \ f_1 \in G(\nu_1), \ \nu_1 \in U\Psi, \nu(a+0)
\in (0,\infty); \ f_2 \in G(\nu_2), \ \nu_2 \in U\Psi, \nu_2(b-0)
\in (0,\infty), $
then $ f_1 \in G(\psi), \ f_2 \in
G(\psi), $ hence $ f \in G(\psi). $ \par
 Inversely, let $ \psi \in U\Psi, \psi(a+0) = \psi(b-0) = 0. $
Let $ p_0 $ be a some number inside the interval $ (a,b) $ such that
$$
\psi(p_0) = \min_{p \in (a,b)} \psi(p) \stackrel{def}{=} C.
$$
Define
$$
\nu_1(p) = \psi(p) \cdot I(p \in (a,p_0)) + C \cdot I(p \in [p_0,b)),
$$
$$
\nu_2(p) = C \cdot I(p \in (a,p_0)) + \psi(p) \cdot I(p \in [p_0,b)).
$$
If $ f \in G(\psi), $ then
$$
f(x) = f(x) I(|f(x)| \ge 1) + f(x) I(|f(x)| < 1) = f_1 + f_2,
$$
where by virtue of Tchebychev's inequality:
$ \mu \{x: |f(x)| \ge 1 \} < |f|_p < \infty $ for some $ p \in (a,b) $
it follows that
$ f_1 \in G(\nu_1); $ and since $ \forall q > p, A \in \Sigma $
$$
\int_A |f_2|^q \mu(dx) \le \int_A |f_2|^p \mu(dx),
$$
we obtain $ f_2 \in G(\nu_2). $ \par

 It is evident by virtue of Liapunov's inequality that in the bounded
case $ \mu(X) = 1: \ G(\psi) = G(\nu_1). $ \par

 We denote by $  G^o = G^o_X(\psi), \ \psi \in U\Psi $ the closed subspace of
 $ G(\psi), $ consisting on all the functions
 $ f, $ satisfying the following condition:
$$
\lim_{p \to a+0} |f|_p / \psi(p) = \lim_{p \to b-0} |f|_p /\psi(p) = 0,
$$
in the case $ \psi(a+0)=\infty, \ \psi(b-0) = \infty; $
$$
\lim_{p \to b-0} |f|_p/\psi(p) = 0
$$
in the case $ \psi(a+0) < \infty, \ \psi(b-0) = \infty; $
$$
\lim_{p \to a+0} |f|_p/\psi(p) = 0
$$
in the case $ \psi(a+0) = \infty, \ \psi(b-0) < \infty; $ and by
 $ GB =  GB(\psi) $ the closed span in the norm $ G(\psi) $
the set of all the bounded measurable functions with finite support:
$ \mu(supp \ |f|) < \infty. $ \par
 Another definition:  for a two functions $ \nu_1(\cdot), \ \nu_2(\cdot)
\in U\Psi $ we will write $ \nu_1 << \nu_2, $ iff
$$
\lim_{p \to a+0} \nu_1(p) /\nu_2(p) = \lim_{p \to b-0} \nu_1(p)/
\nu_2(p) = 0
$$
in the case $ \nu_2(a+0) = \nu_2(b-0) = \infty $ etc. \par
 If for some $ \nu_1(\cdot), \nu_2(\cdot) \in U\Psi, \ \nu_1 <<
\nu_2 $ and  $||f||G(\nu_1) < \infty, $ then $ f \in G^0(\nu_2).$
 Moreover, if there exists a sequence of functions $ f_n, f_{\infty} $
such that for some $ \nu_1 \in G(\psi, a,b) $
$$
\forall p \in (a,b) \ \Rightarrow |f_n - f_{\infty}|_p \to 0, n \to \infty
$$
and $ \sup_{n \le \infty} ||f_n||G(\nu_2) < \infty,$
then $ ||f_n - f_{\infty}||G(\nu_1) \to 0. $ \par
 We consider now some important examples.
Let $ X = R, \ \mu(dx) = dx, 1 \le a < b < \infty,
\gamma = const > -1/a, \ \nu = const >-1/b, \ p \in (a,b), $
$$
f_{a,\gamma} = f_{a,\gamma}(x) = I(|x| \ge 1) \cdot |x|^{-1/a}
(|\log |x| \ |)^{\gamma},
$$
$$
g_{b,\nu} = g_{b,\nu}(x) = I(|x| < 1) \cdot |x|^{-1/b} |\log x|^{\nu},
$$
$$
h_m(x) = (\log |x|)^{1/m} I(|x|<1), \ m = const > 0,
$$
$$
 f_{a,b;\gamma,\nu}(x) = f_{a,\gamma}(x) + g_{b,\nu}(x), \
 g_{a,\gamma,m}(x) = h_m(x) + f_{a,\gamma}(x),
$$

$$
\psi^p_{a,b;\gamma,\nu}(p) = 2(1-p/b)^{-p\nu-1} \ \Gamma(p \gamma + 1) +
2 (p/a - 1)^{-p \gamma-1} \Gamma(p \nu + 1),
$$
$$
\psi^p_{a,\gamma,m}(x) = 2(p/a-1)^{-p\gamma-1} \Gamma(p \gamma+1) +
2 \Gamma((p/m) + 1),
$$
$ \Gamma(\cdot) $ is usually Gamma - function. \par
 We find by the direct calculation:
$$
\left| f_{a,b;\gamma,\nu} \right|^p_p = \psi^p_{a,b;\gamma,\nu}(p);
\ \left|g_{a,\gamma,m} \right|^p_p = \psi^p_{a,\gamma,m}(p).
$$
 Therefore,
$$
\psi_{a,b;\gamma,\nu}(\cdot) \in \Psi(a,b), \ \psi_{a,\gamma,m}(\cdot) \in
\Psi(a,\infty).
$$
Further,
$$
f_{a,b;\gamma,\nu}(\cdot) \in G(a,b; \gamma+ 1/a, \nu +1/b) \setminus
G^o(a,b; \gamma+ 1/a, \nu + 1/b),
$$
$$
g_{a,\gamma,m}(\cdot) \in G \setminus G^0(a,\infty; \gamma+ 1/a, -1/m),
$$
and $ \forall \Delta \in (0, 1)  \ f_{a,b,\alpha,\beta} \notin $
$$
 G(a,b;(1-\Delta)(\gamma+1/a),
\nu+1/b)) \cup G(a,b;1/a, (1-\Delta)(\nu+1/b),
$$

$$
g_{a,\gamma,m}(\cdot) \in G \setminus G^o(a,\infty; \gamma +
1/a; - 1/m).
$$
 Another examples. Put
$$
f^{(a,b; \alpha,\beta)}(x) = |x|^{-1/b}
\exp \left(C_1|\log x|^{1-\alpha} \right) I(|x| < 1) +
$$
$$
I(|x|\ge 1) |x|^{1/a} \exp \left(C_2 (\log x)^{1 - \beta} \right);
$$
$ 1 \le a < b < \infty; \alpha,\beta = const \in (0,1). $  We have:
$$
\log \left|f^{(a,b;\alpha,\beta)}(\cdot) \right|_p \asymp
(p-a)^{1-1/\alpha} + (b-p)^{1-1/\beta}, \ p \in (a,b).
$$
{\bf Theorem 1.} {\it The spaces } $ G(\psi) $ {\it with respect to the
ordinary operations  and introdused norm } $ ||\cdot||G(\psi) $ {\it
are Banach spaces.} \par
 We need only to prove the completness of $ G(\psi) \ - $ spaces.
Denote
$$
\epsilon(n,m) = ||f_n - f_m||G(\psi), \ \epsilon(n) = \sup_{m \ge n}
\epsilon(m,n),
$$
and assume that $ \lim_{n,m \to \infty} \epsilon(m,n) = 0; $ then
$ \lim_{n \to \infty} \epsilon(n) = 0. $ Let $ p(i), i = 1,2,\ldots $
be the (countable) sequence
 of {\it all } rational numbers of interval $ (a,b). $ We have from the direct
definition of our spaces:
$$
\forall p \in (a,b) \ \Rightarrow
|f_n-f_m|_{p(i)} \le \epsilon(n,m) \psi(p(i)).
$$
 As long as the spaces $ L(p(i)) $ are complete, we conclude that there
exist a functions $ f^{(i)}, \ f^{(i)} \in L(p(i)) $ such that
$$
|f_n - f^{(i)}|_{p(i)} \le \epsilon(n)\psi(p(i)) \to 0, \ n \to \infty.
$$
It is evident that
$$
\mu \{x: \forall i \ f^{(i)}(x) \ne f^{(1)}(x) \} = 0,
$$
i.e. $ f^{(i)}(x) = f^{(1)}(x) \ \mu - $ almost everywhere. Hence
$ \forall i = 1,2,\ldots $
$$
|f_n - f^{(1)}|_{p(i)} \le \epsilon(n)\psi(p(i)),
$$
$$
\forall p \in (a,b) \ \Rightarrow |f_n - f^{(1)}|_{p} \le \epsilon(n)\psi(p),
$$
$$
||f_n - f^{(1)}||G(\psi) =
\sup_{p \in (a,b)} |f_n - f^{(1)}|_p / \psi(p) \le \epsilon(n) \to 0,
$$
$ n \to \infty. $ This completes the proof of theorem 1.\par
 Moreover, the spaces $ G(\cdot) $ are rearrangement invariant (r.i.) spaces
with the fundamental function
$$
\phi(G,\delta) \stackrel{def}{=} ||I(A)||G, \ A \in \Sigma, \ \mu(A)
= \delta \in (0,\infty).
$$
 In our case, for the spaces $ G(\psi), \ \psi(\cdot) \in U \Psi,
\ supp \ \psi = (a,b), \ b \le \infty $ we have:
$$
\phi(G(\psi), \delta) = \sup_{p \in (a,b)} \left[ \delta^{1/p} /\psi(p)
\right].
$$
 Note that in the case $ b < \infty $
$$
\delta \le 1 \ \Rightarrow C_1 \delta^{1/a} \le \phi(G,\delta) \le C_2
\delta^{1/b},
$$
$$
\delta > 1 \ \Rightarrow C_3 \delta^{1/b} \le \phi(G,\delta) \le C_4
\delta^{1/a}.
$$
Moreover, $ \lambda \in (0,1) \ \Rightarrow $
$$
\lambda^{1/b} \phi(G, \delta) \le \phi(G,\lambda \delta)  \le \lambda^{1/a}
\phi(G,\delta);
$$
$$
 \lambda > 1 \ \Rightarrow \lambda^{1/b} \phi(G,\delta) \le
\phi(G,\lambda \delta) \le \lambda^{1/a} \phi(G,\delta).
$$
 For instance, define in the case $ b < \infty \  \delta_1 =
\exp(\alpha h^2/(h-a)), \ \delta \ge \delta_1 \ \Rightarrow $
$$
p_1 =p_1(\delta) = \log \delta/(2 \alpha) -\left[0.25 \alpha^{-2}\log^2
\delta -a \alpha^{-1}\log \delta \right]^{1/2},
$$
$$
\phi_1(\delta) = \delta^{1/p_1}(p_1-a)^{\alpha};
$$
$$
\delta \in (0,\delta_1) \ \Rightarrow
\phi_1(\delta) = \delta^{1/h}(h-a)^{\alpha};
$$
$$
\delta_2 = \exp(-h^2\beta/(b-h)), \ \delta \in (0,\delta_2) \ \Rightarrow
$$
$$
p_2 = p_2(\delta) = - |\log \delta|/2 \beta +
\left[ \log^2 (\delta/(4 \beta^2)) + b |\log \delta|/ \beta \right]^{1/2},
$$
$$
\phi_2(\delta) = \delta^{1/p_2(\delta)} (b-p_2(\delta))^{\beta};
$$
$$
\delta \ge \delta_2 \ \Rightarrow \phi_2(\delta) = \delta^{1/h} (b-h)^{\beta}.
$$
 We obtain after some calculations:

$$
b < \infty \ \Rightarrow \phi(G(a,b;\alpha,\beta),\delta) = \max
\left[ \phi_1(\delta), \phi_2(\delta) \right].
$$
 Note that as $ \delta \to 0+ $
$$
\phi(G(a,b,\alpha,\beta), \delta) \sim
( \beta b^2/e)^{\beta} \ \delta^{1/b} \ |\log \delta|^{-\beta},
$$
and as $ \delta \to \infty  $
$$
\phi(G(a,b,\alpha,\beta),\delta) \sim
(a^2 \alpha/e)^{\alpha} \delta^{1/a} \ (\log \delta) ^{-\alpha}.
$$
 In the case  $ b = \infty, \beta < 0 $ we have denoting
$$
\phi_3(\delta) = ( \beta /e)^{\beta}  \ |\log \delta|^{-|\beta|}, \ \delta \in
(0,\exp(-h |\beta|)),
$$
$$
\phi_3(\delta) = h^{-|\beta|} \delta^{1/h}, \ \delta \ge \exp(-h |\beta|):
$$
$$
\phi(G(a,\infty; \alpha, - \beta), \delta) = \max(\phi_1(\delta), \phi_3
(\delta)),
$$
and we receive as $ \delta \to 0+ $ and as $ \delta \to \infty $
correspondingly:
$$
\phi(G(a,\infty;\alpha, - \beta), \delta) \sim (\beta)^{|\beta|}
|\log \delta|^{-|\beta|},
$$
$$
\phi(G(a,\infty; \alpha,-\beta),\delta) \sim ( a^2 \alpha/e)^{\alpha}
\ \delta^{1/a} (\log \delta)^{-a}.
$$

\vspace{3mm}

{\bf 2 \ Connection with another r.i. spaces.} \\

\vspace{2mm}
{\bf Theorem 2.} {\bf A.}  {\it Let} $ \psi(\cdot) \in EX\Psi, $ {\it
such that }
$ \exists g:X \to R, \ \psi(p) \asymp |g(\cdot)|_p, \ p \in (a,b). $
{\it Denote }
$$
N^{(-1)}(1/\delta) = 1/(\phi(G(\psi), \delta)), \ \delta \in (0,\infty),
$$
{\it where } $ N^{(-1)} $ {\it denotes the left inverse function to the }
$ N(\cdot) $ {\it on the set } $ R_+. $ {\it If }
$$
 \forall \epsilon > 0 \
\int_X N(\epsilon |g(x)|) \ \mu(dx) = \infty, \eqno(2.1)
$$
{\it then the space } $ G(\psi) $ {\it is not equivalent to arbitrary
Orlicz's space } $ Or(X,\mu, \Phi). $ \\
{\bf B.} {\it Denote } $ T(x) = (1/\phi(x))^{(-1)}. $  {\it If }
$$
\sup_{p \in supp \ \psi} \left[ \left(\int_0^{\infty} x^{p-1} T(x) dx \right)
/\psi(p) \right]^{1/p} = \infty, \eqno(2.2)
$$
{\it then the space } $ G(\psi) $ {\it  is not equivalent to arbitrary
Marcinkiewicz's space } $ M(\theta).$ \\
{\bf C.} {\it Let } $ \psi(\cdot) \in U\Psi, \ supp \ \psi = (a,b), \
1 \le a < b < \infty. $ {\it Then the space } $ G(\psi) $ {\it is not
equivalent to arbitrary Lorentz space }   $ L(\chi). $ \\
{\bf Proof. A.} Assume conversely, i.e. that $ G(\psi) \sim Or(\Phi), $ where
$ Or(\Phi) $ is some Orlicz's space on the set $ (X, \Sigma, \mu) $ with
 corresponding (convex, even, $ \Phi(0) = 0 $ etc.)
Orlicz's function $ \Phi(u), u \in R. $  Since for $ A \in \Sigma, \mu(A) \in
(0, \infty) $
$$
\phi(Or(\Phi); \mu(A)) = ||I(A)||Or(\Phi) = 1/ \left[ \Phi^{-1}(1/\mu(A)) \right],
$$
we conclude that $ \Phi(u) =  N(u). $ It is evident that $
g(\cdot) \in G(\psi) = Or(\Phi), $ but $ g(\cdot) \notin Or(\Phi)
$  by virtue of our condition (2.1). This contradiction proves the
assertion {\bf A.} \par
 As a consequence: \\
{\bf Lemma 2.} The space $ G(a,b;\alpha,\beta) $ are equivalent to
the Orlicz's space {\it only in the case } $ \alpha = 0, b =
\infty, \beta < 0. $ \par
 (The case $ \alpha = 0, b = \infty, \beta < 0 $ was considered in
[12].) \par
{\bf Proof B.} Assume conversely, i.e. that the space $ G(\psi) = G(\psi,a,b) $
is equivalent
to some  Marcinkiewicz space $ M(\theta) $ over the our measurable space
$ (X,\mu). $  Recall here that in the considered case $ a \ge 1; b > a $
the norm of a function
$ f: X \to R $ in the Marcinkiewicz space may be calculated by the formula
$$
||f||M(\theta) =
\sup_{\delta > 0} \left[\theta(\delta) T^{(-1)}(f,\delta) \right]
$$
and that the fundamental function for the $ M(\theta) $ space in
equal to
$$
\phi(M(\theta), \delta)= 1/\theta(\delta),
$$
(see, for example, [21], p. 187). Therefore, if the space $ G(\psi) $ is
equivalent to some Marcinkiewicz space $ M(\theta), $ then
$$
\theta(\delta) = \delta/\phi(G(\psi),\delta)).
$$
 Let us consider the function $ f:X \to R $ with the tail - function
$ T(f,x) = T(x), $ then $ f \in M(\theta), $ but it follows from our
condition (2.2) that $ f \notin G(\psi). $  \par
 For example, all the spaces $ G(a,b;\alpha,\beta) $ are not equivalent to
arbitrary Marcinkiewicz space.\par
{\bf Proof C } is very simple, again by means of the method of "reduction in
absurdum". Suppose $ G(\psi) \sim L(\chi), $ where $ L(\chi) $ denotes the
Lorentz space with some (quasi) - concave  generating
function $ \chi(\cdot). $ Since
$$
\phi(L(\chi),\delta) = \chi(\delta) \to 0, \delta \to 0+
$$
and $ \chi(\delta) \to \infty, \ \delta \to \infty, $ we conclude that
the space $ L(\chi) = G(\psi) $ is separable ([22], p. 150.) But we will
prove further (in the section 4)
that the space $ G(\psi) $ are non - separable.\par

\vspace{3mm}

{\bf 3 \ Norm's absolute continuity.} \\

\vspace{3mm}

 We will say that the function
$ f \in G(\psi), \ \psi \in U\Psi $ has {\it absolute continue norm}
and write $ f \in GA(\psi), $ if
$$
\lim_{\delta \to 0} \sup_{A: \mu(A) \le \delta} ||f \ I_A||G(\psi)
=0.
$$
The subspaces $ GA(\psi), GB(\psi), G^0(\psi) $ are closed subspaces of
space $ G(\psi). $ \par
{\bf Theorem 3.}  {\it Let } $ \psi \in U\Psi. $ {\it Then }
$$
G^0(\psi) = GB(\psi) = GA(\psi).
$$
 For example, if $ \min(\alpha, |\beta|) > 0, 1 \le a < b \le \infty, $ then
$$
G^o(a,b; \alpha,\beta) = GB(a,b;\alpha,\beta) = GA(a,b;\alpha,\beta).
$$
{\bf Proof.} The inclusions $ GB \subset GA, GA \subset G^o $ are obvious.
Let now $ f \in G^0; $ for simplicity we will suppose
 $ b < \infty, \mu(X) = 1. $ Then $ \lim_{p \to b-0} |f|_p /
\psi(p) = 0.$ Let $ \epsilon > 0. $ We have:
$ ||f \ I(|f| \ge N)||G \le $
$$
\sup_{p \in [1,b-\delta]}|f \ I(|f| \ge N|_p/\psi(p) +
 \sup_{p \in (b-\delta,b)} |f|_p /\psi(p) = \Sigma_1 + \Sigma_2;
$$

$$
\Sigma_2 \le \sup_{p \in [b-\delta,b)} |f|_p /\psi(p) \le \epsilon/2
$$
for some $ \delta \in (0,b) $ by virtue of condition  $ f \in G^o. $

 Further, there exists a value $ N \ge 1 $ such that
$$
\Sigma_1 \le C |f \ I(|f| \ge N)|_{b-\delta} \le \epsilon/2
$$
as long as $ f \in L_{b-\delta}. $ Following,
$ f \in GB; $ thus $ G^0 \subset GB. $ \par
 Now we prove the inverse embedding. Let
$ f \in GB, \epsilon > 0. $ Then $ \exists g, \ \sup_x |g(x)| = B < \infty,
\forall p \in [1,b) \ \Rightarrow |f-g|_p /\psi(p) < \epsilon/2, $
$$
|f|_p \le |g|_p + 0.5 \epsilon \psi(p), \ p \in [1,b);
$$
$$
|f|_p/\psi(p) \le |g|_p /\psi(p) + 0.5 \epsilon < 0.5 \epsilon
+ 0.5 \epsilon \le \epsilon, \ |p - b| < \delta
$$
for sufficiently small value $ \delta.$  Theorem 3 is proved. \par
 We investigate here the {\it sufficient} condition for the convergence

$$
||f_n - f_{\infty}||G(\psi) \to 0, \ n \to \infty.  \eqno(3.1)
$$
Assume at first that (the necessary condition)
$$
{\bf A}. \forall p \in (a,b) \ |f_n - f_{\infty}|_p \to 0, \ n \to \infty.
$$
{\bf Theorem 4.} {\it Let } $ f_n, f_{\infty} \in G(\psi). $
 {\it Assume that (in addition to the condition } {\bf A}) \\
{\bf B.} $ \exists \psi_2(\cdot) \in U\Psi, \ \psi << \psi_2, $
such that
$$
\sup_{n \le \infty} ||f_n||G(\psi_2) < \infty.
$$
{\it Then the convergence (3.1) holds. }  \\
{\bf Proof.} We need use the following auxiliary well - known facts.\\
{\bf 1.} Let $ 1 \le a < b \in (1, \infty). $  We assert that
$$
 \sup_{p \in (a,b)} |f|_p < \infty \ \Leftrightarrow
 \max(|f|_a, |f|_b) < \infty.
$$
This proposition follows from the formula
$$
|f|_p^p = p \ \int_0^{\infty} z^{p-1} T(f,z) dz,
$$
Tchebychev's inequality and Fatou's lemma.\\
{\bf 2.}  Let $ 1 \le p(1) \le p \le p(2) < \infty, \max(|f|_{p(1)},|f|_{p(2)})
< \infty. $ Then $ |f|_p \le $
$$
 |f|_{p(1)}^{(p(2)-p)/(p(2)-p(1))} \cdot |f|_{p(2)}^{(p-p(1))/
(p(2)-p(1))} \stackrel{def}{=} Z(p,p(1),p(2); |f|_{p(1)},|f|_{p(2)}).
$$
 Proposition {\bf 2} follows from  H\"older's  inequality.\par
 It is sufficient to investigate the case
$ b < \infty; $
another cases may be proved analogously. Consider the norm
$$
\Sigma \stackrel{def}{=} ||f_n - f_{\infty}||G(\psi) =
\sup_{p \in (a,b)} |f_n - f_{\infty}|_p /\psi(p).
$$
Let $ \epsilon = const > 0. $
 We have: $ \Sigma \le \Sigma_1 + \Sigma_2 + \Sigma_3, $ where
$ \Sigma_1 = $
$$
\sup_{p \in (a,a+ \delta)} |f_n-f_{\infty}|_p /\psi(p) \le
$$

$$
\sup_p \left[|f_n - f_{\infty}|_p/\psi_2(p) \right] \cdot
\sup_{p \in (a,a+\delta)} \psi(p)/\psi_2(p) \le C(a,\delta)  < \epsilon/3,
$$
if $ \delta = \delta(\epsilon) $ is sufficiently small. Further,
$ \Sigma_3 =  $
$$
\sup_{p \in (b-\delta,\delta)} \left[|f_n - f_{\infty}|_p/ \psi_2(p) \right]
\cdot \sup_{p \in (b-\delta,b)} \left[\psi(p)/\psi_2(\delta) \right] \le
C(b,\delta) < \epsilon/3.
$$
 Finally, $ \Sigma_2 \le $
$$
 \sup_{p \in (a+\delta, b-\delta)} |f|_p /\psi(p) \le C Z
\left(p,a+\delta,b-\delta,|f_n - f_{\infty}|_{a+\delta},
|f_n - f_{\infty}|_{b-\delta} \right)
$$
$ < \epsilon/3 $ for sufficiently large values $ n. $ \par
 Analogously may be proved the following assertion about the $ G(\psi) $
convergence.\\
{\bf Lemma 3.} {\it If the sequence of a functions } $ \{f_n(\cdot) \} $
{\it convergens in all the } $ L_p $ {\it norms: }
$$
\forall p \in (a,b) \ \Rightarrow \lim_{n \to \infty}|f_n - f_{\infty}|_p
= 0
$$
{\it and has a uniform absolute continuous norms in the } $ G(\psi) $
{\it space: }
$$
\lim_{\delta \to 0+} \sup_{n \le \infty} \sup_{A: \mu(A) \le \delta}
||f_n \ I(A)||G(\psi) = 0,
$$
{\it then } $ ||f_n - f_{\infty}|| G(\psi) \to 0, \ n \to \infty. $ \\
{\bf Theorem 5.}  {\it Let } $ \psi \in U\Psi. $ {\it We assert that }
 $ ||f||G/G^o = ||f||G/GA = $
$$
||f||G/GB = \inf_{g \in GB} ||f-g||G = \overline{\lim}_{\delta \to 0+}
\sup_{A: \mu(A) \le \delta} ||f I(A)||G.
$$
{\it Here the notation}  $ G/G^o $ {\it denotes the factor -
space.} \par
{\bf Proof. } Suppose for simplicity  $ b \in
(1,\infty), \  \mu(X) = 1, G = G(\psi), \psi(a+0) <\infty,
\psi(b-0) = \infty; \ f \in G \setminus G^o. $
Put
$$
\gamma = \overline{\lim}_{\delta \to 0} \sup_{A: \mu(A) \le \delta}
||f \ I(A)||G > 0.
$$

 Let also $ g = g(x) $ be a measurable bounded function: $ \sup_x |g(x)| =
B \in (0,\infty); k = const \ge 2. $ We conclude using the elementary
inequality: $ X \ge k Y > 0, k > 2, Y \le B = const \ \Rightarrow $
$$
\frac{(X-Y)^p}{X^p-B^p} \ge \frac{(k-1)^p}{k^p-1}:
$$

$$
||f-g||G \ge \sup_{p \in [1,b)}  \left[ \int_{ \{ x: |f(x)| >
k |g(x)| \}} |f(x) - g(x)|^p \ \mu(dx) \right]^{1/p}/\psi(p) \ge
$$
$$
\overline{\lim}_{p \to b-0}
 \left[ \int_{ \{|f(x)| \ge k B \} } (k-1)^p
(k^p-1)^{-1} (|f|^p - B^p) \ \mu(dx) \right]^{1/p}/\psi(p) \ge
$$
$$
(k-1)(k^b - 1)^{-1/b} \ \overline{\lim}_{\delta \to 0} ||f \ I(A)||G =
(k-1)(k^b -1)^{-1/b} \ \gamma.
$$
 Since the value of $ k $ is arbitrary, it follows from the last inequality
that $ ||f-g||G \ge \gamma; $ this proves that
 $ \inf_{g \in GB} ||f-g||G \ge \gamma; $  the inverse inequality is
evident. \par

\vspace{3mm}

{\bf 4 \ Non - separability.}

\vspace{3mm}

 Recall that $ \psi(a+0) = \infty $ or $ \psi(b-0) = \infty. $ \par
{\bf Theorem 6.} {\it The spaces} $ G(\psi), \ \psi \in U\Psi $
{\it are non - separable.} \par {\bf Proof.} The assertion of
theorem 6 is trivial if the metric space $ (\Sigma, \rho(A,B)), \
\rho(A,B) = \arctan(\mu(A \Delta B)) $ is non - separable.
Therefore by virtue of Rockling's theorem we can suppose the space
$ X $ is equipped by the distance $ d = d(x_1,x_2) $ such that the
space $ (X,d) $ is complete and separable, the measure $ \mu $ is
Borelian and diffuse. \par
 Conversely, assume that the space $ G(\psi) $ is separable. Let
$ \{u_n(x)\} $ be a enumerable dense subset of $ G(\psi). $
 By virtue of Lusin's and Prokhorov's theorems we conclude that
  there exists a compact
subset  $ Y $ of $ X $ with $ \mu(Y) > 0 $ such that on the subspace $ Y $
all the functions $ u_n(x) $ are continuous. We consider now the space
$ G(Y,\psi). $ The functions $ \{u_n(x) \}, x \in Y $ belong to the
space $ G^o_Y(\psi). $  Let $ w(x), x \in Y, $ be some function from
the space  $ G_Y(\psi) \setminus G^o_Y(\psi) $ and define $ w(x) = 0, \
x \in X \setminus Y. $ We get:

$$
\inf_n||w-u_n||G_X \ge \inf_n||w-u_n||G_Y \ge \inf_{g \in GB_Y} ||w-g||G_Y > 0,
$$
in contradiction. This completes the proof of theorem 3. \par
 Our proof of theorem 3 is the same as proof of non - separability of
Orlicz's spaces ([1],  p. 103; [2], p. 127).  \par

\vspace{3mm}

{\bf 5 \ Adjoint spaces.} \\

\vspace{2mm}

 The complete description of spaces conjugated to
$ \cap_p L_p, $ see in [3], [4]. The spaces which are conjugate to Orlicz's
spaces are described in [2], p. 128 - 132. The structure of spaces
$ G^*(\psi) $ is analogous. \par
 It is easy to verify that the structure of linear continuous functionals
over the space $ G^0(\psi) = GA = GB $ is follows: $ \forall l \in
G^{0*}(\psi) \Rightarrow  \exists g: X \to R, $
$$
 l(f) = \int_X f(x) g(x) \ \mu(dx).
$$
 We investigate here only some necessary
conditions for the inclusion $ g \in G^*(\psi). $   Notation:
 $  l_g(f) = \int_X f(x) g(x)  \mu(dx). $  Note at first that if
$  \psi \in U\Psi(a,b),  \ q \in (b/(b-1), a/(a-1)) $ and $ g \in  L_q,$ then
 $ g \in G^*(\psi). $ \\
{\bf Theorem 7.} {\it If} $ g \in G^*, $ {\it then} $ \exists  K = K(g) <
\infty \ \Rightarrow $
$$
\forall z > 0 \ \Rightarrow \int_z^{\infty} T(g,u) du \le K \phi(G, T(g,z)).
$$
{\bf Proof.} Let $ l_g \in G^*. $  It follows from uniform
boundedness principe that $ \forall f \in G \ \Rightarrow $
$$
|l_g(f)| =  \left|\int_X f(x) \ g(x) \mu(dx) \right| \le K ||f||_G.
$$
Put $ f = I_A(x), A \in \Sigma, A = \{x: |g(x)| > z \}, \ z > 0; $ then
$$
\int_z^{\infty} T(g,u) du =
 \int_X |g(x)| I(|g(x)| > z) \ \mu(dx) \le K \phi_G(T(g,z)).
$$

 Let now $ \psi \in U\Psi, \ supp \ \psi = (a,b), b < \infty.$ Introduce the
following $ N \ - $ Orlicz function
$$
N_{\psi}(u) = \sup_{p \in (a,b)} \left[ |u|^p \psi^{-p}(p) \right],
$$
then the following implication holds:
$$
\exists \epsilon > 0 \ \int_X N_{\psi}( \epsilon f) \mu(dx) < \infty \
\Rightarrow f \in G(\psi).
$$
 Therefore, the Orlicz's space $ Or(N,X,\mu) $ is subspace of $ G(\psi). $
Following,
$$
(G(\psi))^{*} \subset (L(N_{\psi}))^*.
$$
 Since the function
$ N_{\psi}(u) $ satisfies the $ \Delta_2 $ condition, the adjoint
space  $ (L(H_{\psi}))^* $ may be described as a new Orlicz's space,
namely
$$
(L(N_{\psi}))^* = L \left(\Phi_{\psi} \right), \ \Phi_{\psi}(u) =
\sup_{z \in R} (uz-N_{\psi}).
$$
Thus, we obtained: $ \psi \in U\Psi(a,b), \ 1 \le a < b < \infty \ \Rightarrow $
$$
(G(\psi))^* \subset L \left(\Phi_{\psi} \right).
$$

\vspace{2mm}

{\bf 6 \ Tail behavior.} \\

\vspace{2mm}
 Let $ f \in G(\psi), \ \psi \in U\Psi(a,b), b \le \infty. $
It follows from Tchebychev's inequality that
$$
T(f,u) \le \inf_{p \in (a,b)} \left[ ||f||^p \psi^p(p)/u^p  \right],  \ u > 0.
$$
 Conversely,
$$
|f|^p_p = p \int_0^{\infty} u^{p-1} T(f,u) du, \ p \ge 1;
$$
 therefore
$$
||f||G(\psi) =  \sup_{p \in supp \ \psi} \left[ p \ \left[
\int_0^{\infty} u^{p-1}
T(f,u) \ du \right]^{1/p} \ /\psi(p) \right].
$$
 In the particular case the spaces $ G(a,b;\alpha,\beta) $ we obtain
after simple calculations: \\
{\bf Theorem 8. A}. {\it Let} $ f \in G(a,b;\alpha,\beta), 1 \le a < b
< \infty. $ {\it Then }
$$
u \in (0, 1/2) \ \Rightarrow T(f,u) \le C_1(a,b,\alpha,\beta)
|\log u|^{ a \alpha} u^{-a}; \eqno(5.1)
$$
$$
u \ge 2 \ \Rightarrow T(f,u) \le C_2(a,b,\alpha,\beta) (\log u)^{b \beta}
u^{-b}. \eqno(5.2)
$$
{\bf B}. {\it Conversely, suppose} $ \exists a,b, 1 \le a < b < \infty, \gamma,
\tau \ge 0, C_j > 0 $ such that
$$
T(f,u) \le C_1|\log u|^{\gamma} u^{-a}, \ u \in (0,1/2); \  T(f,u) \le
C_2 (\log u)^{\tau} u^{-b}, \ u \ge 2.
$$
{\it Then } $ f \in G(a,b; \gamma+1, \tau + 1). $ \par
{\bf C}. {\it Let now } $ f \in G(a,\infty; \alpha, - \beta), \ \beta > 0. $
{\it We propose that }
$$
T(f,u) \le C_1 |\log u|^{a \alpha} \ u^{-a}, \ u \in (0,1/2],
$$
$$
T(f,u) \le C_2 \exp \left(-C_3 u^{1/\beta} \right), u \ge 1/2;
$$
{\bf D}. {\it Conversely, if } $ \exists a \ge 1, \beta > 0,
 \gamma \ge 0, $
$$
T(f,u) \le C_1 |\log u|^{\gamma} \ u ^{-a}, u \in (0,1/2),
a = const > 0, \gamma \ge 0,
$$
$$
T(f,u) \le C_2 \exp \left(-C_3 u^{1/\beta} \right), \beta > 0,
$$
{\it then} $ f \in G(a,\infty; \gamma+1, - \beta). $ \par

 Note in addition that at $ \min(\alpha,\beta) > 0, b < \infty $
$$
T(f,u) \sim C_1|\log u|^{a \alpha} u^{-a}, u \to 0+ \Leftrightarrow
|f|_p \sim C_2( p-a)^{ -\alpha}, p \to a+0;
$$
$$
T(f,u) \sim C_3 |\log u|^{b \beta} u^{-b}, u \to \infty \Leftrightarrow
|f|_p \sim C_4 (b - p)^{-\beta}, p \to b-0
$$
(Richter's theorem).\par
 Despite the well - known Richter's theorem, we can show that both the
inequalities (5.1) and (5.2) are exact. Let us consider the following
examples.\par
 Example 5.1. Let $ \mu(X)=1, $  i.e. $ (X,\Sigma,\mu) $ is a probability
space and let $ \mu $
is diffuse. Consider the (measurable) discrete = valued
 function $ f: X \to R $ such that
$$
\mu \{x: f(x) = \exp(\exp(k)) \} = C \exp(\beta b k - b \exp k), k = 1,2,\ldots;
$$
$$
1/C = \sum_{k=1}^{\infty} \exp(\beta b k - b \exp (k)),
$$
and denote $ \gamma = \beta b, \ a(k) = a(k,\gamma,\epsilon) =
\exp(k \gamma - \epsilon \exp(k)),$

$$
\epsilon = b-p \to 0+, \ k(0) \stackrel{def}{=} [\log (\gamma/\epsilon)],
\ x(k) = \exp(\exp(k)),
$$
here $ [z] $ denotes the integer part of $ z. $ We get:
$$
W(\epsilon) \stackrel{def}{=}C^{-1}|f|_p^p  = \sum_{k=1}^{\infty}
a(k,\gamma,\epsilon) \ge
$$
$$
C_2 a(k(0),\gamma,\epsilon) \ge C_3 (b-p)^{-b \beta},
$$
therefore $ |f|_p \ge C_4 (b-p)^{-\beta}. $ \par
 On the other hand, we have at $ k > k(0) $ and $ k < k(0) $ correspondently
$$
a(k+1)/a(k) < \exp(\gamma(e-2)) < 1, \ a(k-1)/a(k) < \exp(-\gamma/e) < 1,
$$
hence
$$
W(\epsilon) \le C_3 a(k(0),\gamma,\epsilon) \le C \epsilon^{-p \beta},
$$
following $ |f|_p \le C_5 (b-p)^{-\beta}, p \in (1, b). $ Thus
$ f \in G(1,b;0,\beta). $ However,
$$
T(|f|,x(k)) > C \exp(b \beta k - b \exp k) = C (\log x(k))^{b \beta}
x(k)^{-b}.
$$
(we used the discrete analog of saddle - point method).\par

Example 5.2. Let $ X = R^1_+, \mu(dx) = dx, Q(k) = \exp(a\alpha k + a \exp(k)),
a = const \ge 1, S(k) = \sum_{l=1}^k Q(l), \ b \in (a,\infty), $
$$
g(x) = \sum_{k=1}^{\infty} \exp(-\exp(k)) \ I(x \in (S(k-1), S(k)]),
$$
$ u(k) = \exp(-\exp(k)). $
We obtain analogously to the example 5.1:
$$
p \in (a,b)  \Rightarrow |g|_p \asymp (p-a)^{-\alpha},
$$
but
$$
T(g,u(k)) \ge C(a,b, \alpha) \ |\log u(k)|^{a \alpha} u(k)^{-a}.
$$


\newpage

 {\bf 7 \ Fourier's transform.}

\vspace{2mm}

 In this section we investigate the boundedness of certain
Fourier's operators, convergence and divergence Fourier's series
and transforms in $ G(\psi) $ spaces. Let
$ X = [-\pi,\pi] $ or $ X = R = (-\infty,\infty), \ \mu(dx) = dx, \
 c(n) = c(n,f) = $
$$
 \int_{-\pi}^{\pi} \exp(i n x) f(x) dx, n = 0,\pm 1,\pm 2 \ldots; \
 2 \pi s_M[f](x)  =
$$
$$
\sum_{ \{n: |n| \le M \} } c(n) \exp(-inx), \ s^*[f] = \sup_{M \ge 1}
|s_M[f]|,
$$

$$
F[f](x) = \lim_{M \to \infty} \int_{-M}^M \exp(i t x) f(t) dt,
$$
$$
F^*[f](x) = \sup_{M > 0} \int_{-M}^M \exp(itx) f(t) dt,
$$
$$
S_M[f](x) = (2 \pi)^{-1} \int_{-M}^M  \exp(-itx) F[f](t) dt,
$$
$$
S^*[f](x) = \sup_{M > 0} |S_M[f](x)|.
$$

 Recall that if $ f \in L_p(R),  p \in [1,2], $ then operators $ F,F^* $
are well defined; for the values $ p > 1, \ f \in L_p $ are well
defined the operators $ s_M, s^*, S_M, S^*.$  \par
 We introduce also for arbitrary $ \psi(\cdot) \in U\Psi,
\ supp \ \psi \supset (1,2],
\ \psi_1(p) = \psi(p/(p-1)), $ for $ s = const \in (1,\infty),
 \psi(\cdot) \in U\psi, \ supp \ \psi \supset (1,s) $
$$
\psi_{(s)}(p) = \psi(sp/(s-p)); \ p = \infty \ \Rightarrow p/(p-1)=+\infty;
$$
for $ \psi \in U\Psi, \ supp \ \psi \supset [1,s/(s-1)), $
$$
\psi^{(s)}(p) = \psi[ps/(s-1)/(p+s/(s-1))].
$$
 Let $ \lambda,\gamma = const \ge 0; $ we denote for $ \psi \in
U\Psi(1,\infty) $
$$
\psi_{\lambda,\gamma}(p) = p^{\lambda + \gamma}\psi(p) \ (p-1)^{-\gamma}.
$$
 It is easy to verify that if $ \psi \in EX\Psi, $ then $
\psi_{\lambda,\gamma} \in EX\Psi. $ \par

 Let $ X,Y $ be a two Banach spaces and let $ F: X \to Y $ be a operator
(not necessary linear or sublinear) defined on the space $ X $ with
values in $ Y. $ The operator $ F $ is said to be bounded from the space
$ X $ into the space $ Y, $  notation:
$$
||F||[X \to Y] < \infty,
$$
if for arbitrary $ f \in X \ \Rightarrow \ ||F[f]||Y \le C \cdot ||f||X. $ \\
{\bf Theorem 9.} {\it Let } $ \psi \in U\Psi,
 (1,2] \subset \ supp \ \psi. $
{\it The operator } $ F $  {\it is bounded from the space} $ G(\psi) $
{\it into the space } $ G(\psi_1): $
$$
||F||[G(\psi) \to G(\psi_1)] < \infty.
$$
{\bf Proof.} We will use the classical result of Hardy - Littlewood -
Young:
$$
|F[f]|_{p/(p-1)} \le C |f|_p, \ p \in (1,2].
$$
 Here $ C $ is an absolute constant. \par
 If $ f \in G(\psi), $ then $ |f|_p \le ||f||G \cdot \psi(p), $ therefore

$$
|F[f]|_p \le \psi(p/(p-1)) \ ||f||G(\psi) = \psi_1(p) \ ||f||G(\psi).
$$
 {\bf Theorem 10. } {\it Let } $ X = [-\pi,\pi],   \psi \in U\Psi,
supp \ \psi \supset (1,\infty). $ {\it  We assert that }
$$
\sup_{M \ge 1} ||s_M||[G(\psi) \to G(\psi_{1,1}) ] < \infty.
$$
 {\bf Proof.} Now we use the well - known result of M.Riesz:
$$
||s_M[f]||[L_p \to L_p] \le C p^2/(p-1), \ p \in (1,\infty).
$$
with absolute constant C. If $ f \in G(\psi), $ then $ |f|_p \le $
$$
 \psi(p) ||f||G(\psi), \ |s_M|_p \le C p^2 |f|G(\psi) /(p-1) =
C ||f||G(\psi) \cdot \psi_{1,1}(p).
$$
 {\bf Corollary 3.} Assume in addition to the conditions of theorem 10 that
$ supp \ \psi \subset (a,b) $ for some $ a = const > 1, a < b = const
< \infty. $ Then
$$
\psi_{1,1} (p) \asymp \psi(p), \ p \in (a,b).
$$
Therefore, in this case
$$
\sup_{M \ge 1} ||s_M||[G(\psi) \to G(\psi) ] < \infty.
$$
{\it However, this assertion does not means that } $ \forall f \in G(\psi) \
\Rightarrow $
$$
\lim_{M \to \infty} ||s_M[f] - f||G(\psi) = 0;
$$
see counterexamples further. If $ \nu(\cdot) \in U\Psi, \ \nu <<
\psi_{1,1}, f \in G(\psi), $ then
$$
\lim_{M \to \infty} ||s_M[f]-f||G(\nu) = 0,
$$
i.e. the sequence $ s_M[f] $ convergent to the function $ f $ in
the $ G(\nu) $ sense.\par
 At the same assertion is true if $ f \in G^0(\psi). $ \\
 The assertion analogous to the assertion of theorem 10 is true for the maximal
Fourier's operator $ s^*, $ Fourier transform $ S_M $ and maximal Fourier
transform $ S^* $ etc. \par

 Namely, in  [13], p. 163  is proved that $ \forall f \in L_p, p
\in (1,2] \ |F^*[f]|_p \le C p^4(p-1)^{-2} |f|_p. $ Following,
$$
||F^*|| [G(\psi) \to G \left(\psi_{2,2} \right)] < \infty.
$$

  Let us show the exactness of theorem 9. Let $ f(x) =
 f_{a,b}(x) = |x|^{-1/b}, |x| \in (0,1); f(x) =
|x|^{-1/a}, |x| \ge 1; G = G(a,b;1/a, 1/b), G^/ = G(b/(b-1),a/(a-1), (b-1)/b,
(a-1)/a);  $  then  $ f \in G. $
 It is easy to calculate that
$ F[f_{a,b}](t) \asymp f_{b/(b-1),a/(a-1)}(t), \ t \in R, $ so
$$
F[f_{a,b}] \in G^/ \setminus G^{/0}.
$$
 This example is true even in the case $ a = 1; $ then $ a/(a-1) + \infty. $ \par
 For the Fourier series  $ \sum_n c(n) \exp(inx) $  it is well known (on the
basis of Riesz's theorem) that
$$
f \in L_p[-\pi,\pi], \exists p > 1 \ \Rightarrow \lim_{M \to \infty} |s_M[f] - f|_p = 0.
$$
 This fact is true also in the  Orlicz's spaces with $ N - $ function
satisfying the
so - called  $ \Delta_2 \cap \nabla_2 $  conditions ([6], p. 196 - 197).
 Conversely, in the exponential Orlicz's spaces there exist a functions $ f,$
belonging to this spaces but such that Fourier series (or integrals) does not
convergent to $ f $ in the Orlicz's norm sense
 [5].  Analogously, this effect appears in $ G(\psi) $
spaces. \\
{\bf Lemma 4. } {\it Let } $ \psi \in EX\Psi, X = [-\pi,\pi]. $
{\it There exists a function }
$ f \in G(\psi) $ {\it for which the Fourier series does not convergence in }
$ G(\psi) $ {\it norm to the function } $ f. $ \\
{\bf Proof. } Since $ \psi \in EX\Psi, $ there exists a function $ f:
X \to R $ for which $ |f|_p \asymp \psi(p), \ p \in (a,b); $   then
$ f \in G \setminus G^0(\psi). $
Assume conversely, i.e.
$$
\lim_{M \to \infty} ||s_M[f] - f||G(\psi) = 0.
$$
 Since the trigonometrical system is bounded, this means that $ f \in G^0, $
in contradiction.\\

\vspace{2mm}

{\bf 8 \ Martingales. } \\

\vspace{2mm}

 Let $ (f_n, F_n) $ be a martingale, i.e. a
monotonically  non - decreasing sequence
of $ F_n \ - $ sigma - subalgebras $ \Sigma $ and  $ F_n $ measurable
functions $ f_n $ such that $ {\bf E} f_{n+1}/F_n = f_n. $ \par
 In this section we will use the probabilistic notations
$$
{\bf E} f = \int_X f(x) \mu(dx), \ |f|_p = {\bf E}^{1/p}|f|^p
$$
and notation $ {\bf E } f /F $  for the conditional expectation. \par
 The $ L_p \ - $  theory of conditional expectations  and theory of
martingales in the case
$ \mu(X) = \infty $ and some applications see, for example, in the
book  [7],  pp. 330 - 347. \par
 The Orlicz's norm estimates for martingales are used in moderne
non - parametrical statistics, for example, in the so - called regression
problem ([10], [11], [12]) etc. Namely, let us consider the following problem.
Given: the observation of a view
$$
\xi(i) = g(z(i)) + \epsilon(i), \ i = 1,2,\ldots,
$$
where $ g(\cdot) $ is inknown estimated function, $ \{\epsilon(i) \} $ is the
errors of measurements and may be an independent random variables or martingal
differences,  $ \{z(i)\} $ is some dense
set in a metric space $ (Z,\rho) $ with Borel measure $ \nu: \ \ z(i) \in Z. $ \par
 Let $ \{\phi_k(z) \} $ be some complete orthonormal sequence of functions,
for example, classical trigonometrical sequence, Legengre or Hermit polynomials
etc. Put
$$
c_k(n) = n^{-1} \sum_{i=1}^n \phi_k(z(i)), \ \tau(N) = \tau(N,n) =
\sum_{k = N+1}^{2N} (c_k(n))^2,
$$
$$
M = argmin_{n \in [1,n/3]} \tau(N), \ f_n(z) = \sum_{k=1}^M c_k(n) \phi_k(z).
$$
 By the investigation of confidence interval for $ ||f_n - f|| $ are used the
exponential bounds for polynomial martingales. \par
 The next facts about martingales in the unbounded case
$ \mu(X) = \infty $  either there are in [7], p. 347 - 351, or are simple
generalization of classical results in the case $ \mu(X)=1 $ ([8], [9]). \\
{\bf 1.} Let the martingale $ (f_n,F_n) $  be a non - negative,
$ c, d = const,
 0 < c < d <  \infty $ and let for some  $ p \ge 1 \ \sup_n |f_n|_p < \infty. $
Denote by $ \nu = \nu(c,d) $  the number of upcrossing of interval
$ (c,d) $ by the (random) sequence $ \{f_n \}. $ Then

$$
{\bf E} \nu \le (d-c)^{-p} \left[ 2^{p-1} \sup_n |f_n|^p_p  + 2^{p-1} c^p
+ (d-c)^p \right].
$$
{\bf 2.} Almost everywhere convergence. If for some
 $  p \ge 1 \ \sup_n|f_n|_p < \infty, $
then  $ \exists f_{\infty}(x) = \lim_{n \to \infty}
f_n(x) \ (mod \ \mu), \ |f_{\infty}|_p < \infty. $ \\
{\bf 3. } Convergence in $ L_p $ norms. If $ \exists p > 1 \ \Rightarrow
\sup_n |f_n|_p < \infty, $ then
$$
\lim_{n \to \infty} |f_n - f_{\infty}|_p = 0.
$$
{\bf 4.} Doob's inequality: $ p > 1 \ \Rightarrow $
$$
p > 1 \ \Rightarrow \left| \sup_n f_n \  \right|_p \le  \sup_n
\left[|f_n|_p \right] \ p/(p-1).
$$
 In the bounded case $ \mu(X) = 1 $ the convergence of martingale
 $ (mod \ \mu) $ is true under (sufficient) condition
$ \sup_n |f_n|_1 < \infty; $  let us show here that in  unbounded
case $ (\mu(X) = \infty) $ our condition is unimproved. Namely, we
consider the sequence of independent identically distributed
functions $ h_j = h_j(x) $ such that for some $ p \ge 1 $
$$
|h_j|_p < \infty;  \ \forall s \ne p, s \ge 1 \ \Rightarrow |h_j|_s = \infty.
$$
Put
$$
f_n(x) = \sum_{j=1}^n 2^{-j} h_j(x), \ F_n = \sigma
\{h_j, j \le n \};
$$
then the convergence $ f_n(\cdot) \ (mod \ \mu) $  is true, despite
$ \forall \ s \ne p \ |f_n|_s = \infty. $ \par
It is proved in the book [10], p. 252, see also [11] that if in some
Orlicz's space $ Or(X,\Sigma, \mu; N) = Or(N), $ with $ \mu(X) = 1 $ and
with the
$ N \ - $ Orlicz's function satisfying $ \Delta_2 \cap \nabla_2 $
condition the martingale $ \{ f_n \} $ is bounded:
$$
\sup_n ||f_n||Or(N) < \infty,
$$
then the martingale  $ \{ f_n \} $ convergent in the correspondent
Orlicz's norm:
$$
\lim_{n \to \infty} ||f_n - f_{\infty}||Or(N) = 0.
$$
 In the article [12]  is showed that in the {\it exponential} Orlicz's spaces
$ Or(N) $
the $ Or(N) $ bounded martingale may divergent. Let us prove
that in the $ Or(N) $
spaces is the same case.\\
{\bf Lemma 5.} {\it Let } $ \psi \in EX\Psi, $  {\it so that} $ \psi(p)
\asymp |f|_p, $ {\it and let the } $ \sigma \ - $ {\it algebra }
$ \sigma(f) $ {\it be an union of finite } $ \sigma \ - $ {\it algebras:}
$$
 \sigma(f) = \cup_{n=1}^{\infty} \sigma_n, \ card(\sigma_n) < \infty
$$
{\it with finite subsets: }
$$
\forall A \in \sigma_n, A \ne X \ \Rightarrow \mu(A) < \infty.
$$
{\it Then there exists a bounded but divergent
in } $ G(\psi) \ - $ {\it sense martingale }
$$
(f_n, F_n): \ \sup_n ||f_n||G(\psi) < \infty, \
\overline{\lim}_{n \to \infty} ||f_n - f_{\infty}||G(\psi) > 0.
$$

{\bf Proof.} Let us consider some function $ f \in G(\psi)
\setminus G^0(\psi).$
 Put $ F_n = \sigma_n, \ f_n = {\bf E} f/F_n; $ then
$ (f_n,F_n) $ is a (regular) bounded martingale:

$$
 \sup_n||f_n||G = \sup_{p \in (a,b)} |f_n|_p /\psi(p) \le \sup_{p \in (a,b)}
|f|_p /\psi(p) = ||f||G < \infty;
$$
we used the Iensen inequality $ |f_n|_p \le |f|_p. $ \par
 Since the sigma - algebras $ \sigma_n $ are finite, $ f_n \in G^0(\psi). $
Suppose $ ||f_n - f||G \to 0, \ n \to \infty, $ then $ f \in G^0, $
in contradiction with choosing $ f. $ \\
{\bf Theorem 11.}  {\it Let } $ (f_n,F_n) $ {\it be a martingale, }
$ \psi \in U\Psi, $
$$
\sup_n ||f_n||G(\psi) < \infty.
$$
{\it Then}
$$
{\bf A.} \ || \sup_n f_n||G \left(\psi_{0,1} \right) < \infty.
$$
{\it Assume in addition that } $ supp \ \psi = (a,b), 1 < a < b \le \infty.$
{\it Then } $ \forall \nu \in U(\psi), \ \nu << \psi_{0,1} $
$$
{\bf B.} \ \lim_{n \to \infty} || f_n - f_{\infty}||G(\nu) = 0.
$$
 {\bf Proof } use the Doob's inequality and
 is the same as in theorem 8 and may be omitted.\par
 For example, let
 $ (f_n,F_n) $  be a martingale, $ 1 \le a < b \le \infty, \
 \sup_n ||f_n||G(a,b;\alpha,\beta) < \infty. $  Then in the case
$ a > 1 $ is true the following implication
$$
|| \sup_n |f_n| \ ||G(a,b; \alpha,\beta) < \infty; \ \forall \Delta > 0
\ \Rightarrow
$$
$$
 \lim_{n \to \infty}||f_n - f_{\infty} ||G(a,b; \alpha +
\Delta, \beta + \Delta [I(b < \infty) - I(b = \infty)]) = 0;
$$
if $ a = 1, $ then
$$
|| \sup_n |f_n| \ ||G(1,b; \alpha+1, \beta) < \infty; \ \forall \Delta > 0 \
\Rightarrow
$$
$$
 \lim_{n \to \infty} ||f_n - f_{\infty}||
G(1,b; \alpha+1 + \Delta,\beta + \Delta [I(b < \infty) - I(b = \infty)]) = 0.
$$
 It is clear that the convergence
 $ f_n \to f_{\infty} $  in the norm $ G(a,b; \alpha,\beta) $ is true also
in the case $ f_{\infty} \in G^o(a,b;\alpha,\beta). $ \par

\vspace{2mm}

{\bf 9 \ Operators.} \\

\vspace{2mm}

 In this section we assume that there is a
measurable space $ (X,\Sigma,\mu) $ and $ Q \ $
is an operator not necessary linear or sublinear defined on the set
  $ \cap_{p \in (a,b)} L_p(X,\mu), 1 \le a = const < b = const \le \infty $
and taking values in the set $ \cap_{p \in (c,d)} L_p(X,\mu). $ We
will investigate the problem of boundedness of operator $ Q $ from
some space $ G(X, \psi) $ into some {\it another } space $
G(X,\nu). $ \par
  The case of Orlicz spaces and certain singular operators was consider in
many publications; see, for example, [18], [19], [20]. \\
 At first we consider the regular operators. \\
{\bf 1.} Define a multiplicative operator
$$
Q_f[g](x) = f(x) \cdot g(x).
$$
 Assume that $ f \in L_s $ for some $ s = const > 1 $ and denote $ t = t(s) =
s/(s-1). $ As long as
$$
|Q_f[g]|_r \le |f|_s \cdot |g|_{rt/(r+t)}, \ r < s,
$$
we conclude: if $ supp \ \psi \supset (t(s), \infty), $ then
$$
||Q_f||[G(\psi) \to G(\psi_{(s)}] < |f|_s, \ \psi_{(s)}(p) = \psi(ps/(s-p)).
$$
{\bf 2.} We consider now the convolution operator (again regular)
$$
Con_f[g](x) = f*g(x) = \int_X g \left(xy^{-1} \right) f(y) \ \mu(dy),
$$
where $ X $ is unimodular Lie's group, $ \mu $ is its Haar measure. Assume
that $ f \in L_s(X,\mu) $ for some $ s = const > 1. $ Using the classical
Young inequality
$$
|f*g|_r \le C(r,s)  |f|_s \cdot  |g|_{rt(s)/(r+t(s))}, \ r > s, C(r,s) < 1,
$$
we observe that
$$
||Con_f|| \left[G(\psi) \to \left(\psi^{(s)} \right) \right] \le |f|_s.
$$
 For example,  if $ \min(\alpha,\beta) > 0, $ then
$$
||Con_f||[G(1,\infty;\alpha, -\beta) \to G(s,\infty; \alpha,0)] \le
C(\alpha, \beta,s)  |f|_s, \ s > 1.
$$
{\bf 3. } Finally we consider some classical singular operators.
Assume that the operator $ Q $ satisfies the following condition:
for some $ \lambda, \gamma = const \ge 0 $ and $ \forall \ p \in (1,\infty) $

$$
|Q[f]|_p \le C \ |f|_p \ p^{\lambda + \gamma} (p-1)^{-\gamma}.  \eqno(8.1)
$$
 There are many singular operators satisfying this condition, for instance,
Hilbert's operator: $ X = (-\pi,\pi) $ (or, analogously, $ X = R), $
$$
H[f](x) = \lim_{\epsilon \to 0+} H_{\epsilon}[f](x),
$$

$$
H_{\epsilon}[f](x) = (2\pi)^{-1} \int_{\epsilon \le |y| \le \pi}
[f(x-y)/\tan(y/2)] dy, \ \lambda = \gamma = 1;
$$
maximal Hilbert's operator
$$
H^*[f](x) = \sup_{\epsilon \in (0,1)} |H_{\epsilon}[f](x)|, \ \lambda =
1, \gamma = 2;
$$
 operators of Caldron - Zygmund: $ \lambda = \gamma = 1, $ Karlesson - Hunt:
$ s^*, S^*; \lambda = 1, \gamma = 3; $ maximal, in particular, maximal
 Fourier, operators, for example,
$$
Q[f](x) \stackrel{def}{=} \sup_{M > 0} \left| \int_R f(t) [\sin(M(x-t))/
(x-t)]  \ dt \right|: \ \lambda = \gamma = 2;
$$
 pseudodifferential operators ([15], p. 143): $ \lambda = 1 = \gamma, $
 oscillating operators ([14], p. 379 - 381) etc.\par
 The following result is obvious.\\
{\bf Theorem  12.} {\it Let } $ \psi \in U\Psi, supp \ \psi = (1,\infty). $
{\it Assume that the operator } $ Q $ {\it satisfies the condition (8.1). }
{\it Then }
$$
||Q|| \left[G(\psi) \to G \left(\psi_{\lambda,\gamma} \right) \right] < \infty.
$$
 Let us consider examples. Assume again that the operator $ Q $ satisfies the
condition (8.1). Then $ Q $ is bounded as operator from the space
$ G(a,b;\alpha,\beta) $ into the space $ G(a,b; \alpha_1,\beta_1),$
where at $ 1 < a < b < \infty \ \Rightarrow \alpha_1 = \alpha,\beta_1 = \beta; $
in the case  $ a = 1, b < \infty \ \Rightarrow \alpha_1 = \alpha + \gamma, \beta_1 =
\beta; $ if $ a > 1, b = \infty $ then $ \alpha_1 = \alpha,
\beta_1 = \beta + \lambda; $  ultimately, for
 $ a=1, b = \infty $ we obtain:
$ \alpha_1 = \alpha + \gamma, \beta_1 = \beta + \lambda. $ \par
 Now we show the exactness of estimations of theorem 12. Consider at first
the singular Hilbert operator for the functions defined on the set
$ (-\pi,\pi). $ \par
 Put now
$$
f(x) = f_d(x) = \sum_{n=2}^{\infty} n^{-1} \log^d n \ \sin(nx), \ d \ge 0.
$$
then  (see [16], p. 184; [17], p. 116])
$ |f(x)| \asymp ( 2 + |\log(|x|)|)^d, \ |f|_p \asymp p^d, p \in [1,\infty),
x \in [-\pi,\pi] \setminus \{ 0 \}; $
$$
\ C H[f](x) = \sum_{n=2}^{\infty} n^{-1} \log^d n \ \cos(nx),
$$
$$
H[f](x) \asymp ( 2 + |\log(|x|)|)^{d+1}, \ |H[f]|_p \asymp p^{d+1}.
$$
 Considering the examples $ d \in (0,1), g = g_d(x) = $
$$
 \sum_{n=1}^{\infty} n^{d-1} \sin(nx), \
C H[g] = \sum_{n=1}^{\infty} n^{d-1} \cos(nx),
$$
we can see that
 $ |g(x)| \asymp |H[g]|(x), x \in R \setminus \{ 0 \}, $ and following
 $ |g|_p \asymp |H[g]|_p, \ p \in (1,\infty). $ \par
 We can built more general examples considering the functions of a view
$$
f(x) = \sum_{n=2}^{\infty} n^{d-1} \ L(n) \ \sin(nx),
$$
 where $ L(n) $ is some slowly varying as $ n \to \infty $ function.
See [17], p. 187 - 188.\par
 The case of Hilbert's transform on the real axis is investigated
analogously. Namely, consider the functions

$$
f(x) = \int_3^{\infty} t^{d-1} \sin(tx) \ dt, \ d \in (0,1),
$$
then (see [17], p.117)  $ C H[f](x) = $
$$
 \int_3^{\infty} t^{d-1} \cos(tx) \ dt, \ | H[f](x)| \asymp
|f(x)| \asymp f_{1/d,1}(x),
$$
following,
$$
H[f](\cdot), \ f(\cdot) \in G \setminus G^o(1,1/d;1,d).
$$

Analogously, considering the example

$$
 f(x) = \int_3^{\infty} t^{-1} \ \sin(tx) dx, \ |f(x)|  \asymp f_{\infty,1}(x),
$$
$ x \in R \setminus \{ 0 \}, $  we observe that $ |H[f](x)| \asymp |\log |x||,
\ |x| \le 1/2; $
$$
f(\cdot) \in G \setminus G^o(1,\infty;1,0), \ |C H[f](x)| \sim |\log |x||, \
x \to 0;
$$
$$
 | H[f](x)| \asymp |x|^{-1}, \ |x| \ge 1/2,
$$
so that $ H[f](\cdot) \in G \setminus G^o(1,\infty;1,1), $ \par
{\bf Acknowledgement.} We are very grateful for support and
attention to prof. V.Fonf, L.Beresansky, and M.Lin (Beer - Sheva,
Israel). This investigation was partially supported by ISF (Israel
Science Foundation), grant $ N^o $ 139/03.\par

\newpage
\vspace{3mm}
{\bf References.} \\
\vspace{2mm}

1. Krasnoselsky M.A., Rutisky Ya.B. Convex functions and Orlicz's
Spaces. P. Noordhoff LTD, The Netherland, 1961, Groningen. 2. Rao
M.M, Ren Z.D. Theory of Orlicz Spaces. - New York, Basel. Marcel
Decker,
 1991. - 449 p. \\
3. Davis H.W., Murray F.J., Weber J.K. Families of $ L_p - $ spaces with
inductive and projective topologies.// Pacific J.Math. - 1970 - v. 34,
p. 619 - 638.\\
4. Steigenwalt M.S. and While A.J. Some function spaces related to
$ L_p. $ //
Proc. London Math. Soc. - 1971. - 22, p. 137 - 163.\\
5. Ostrovsky E., Sirota L. Fourier Transforms in Exponential Rearrangement
Invariant Spaces.// Electronic Publ., arXiv:Math., FA/040639, v.1, - 20.6.2004.\\
6. Rao M.M., Ren Z.D. Application of Orlicz Spaces. -  New York, Basel.
 Marcel Decker, 2002. - 475 p. \\
7. Rao M.M. Measure Theory and Integration. - New York, Basel,
 Marcel Decker, second edition, - 2004. - 781 p. \\
8.  Hall P., Heyde C.C. Martingale Limit Theorems and its Applications. -
USA, New York, Academic Press Inc., 1980, - 473 p.\\
9. Doob J.L. Stochastic Processes. -  New York, John Wiley, - 1953, 671 p.\\
10. Neveu J. Discrete - Parameter Martingales. - Amsterdam - Oxford - New York,
North - Holland Publ. Comp., 1975. - 385 p.\\
11. Peshkir G. Maximal Inequalities of Kahane - Khintchin type in Orlicz
Spaces. // Preprint Series 33, Inst. of Math., 1992, University of Aarhus
(Denmark). - 1992. - 1 - 44 p. \\
12. Ostrovsky E.  Exponential Orlicz Spaces: new Norms and Applications.//
Electronic Publ., arXiv/FA/0406534, v. 1, - 25.06.2004.\\
13. Juan Arias de Reyna.  Pointwiese Convergence Fourier Series. -
New York,
Lect. Notes in Math. 2004. - 286 p. \\
14. Stein E.M. Singular Integrals and Differentiability Properties of
Functions. - Princeton, New York. University Press. 1970. - 572 p. \\
15. Taylor M. Pseudodifferential Operators. - Princeton, New
Jersey,
University Press. 1981. - 473 p. \\
16. Zygmund A. Trigonometrical Series. Volume 1. Cambridge University Press,
1959.\\
17. Seneta E. Regularly Varying Functions. 1985, Moscow edition.\\
18. Bongioanni B., Forzani L., Harboure E. Weak type and restricted weak
type (p,p) operators in Orlicz spaces. (2002/2003). Real Analysis Exchange,
Vol. 28(2), pp. 381 - 394.\\
19. Harboure E., Salinas O. and Viviani B. Orlicz Roundedness for
Certain Classical Operators. Colloquium Matematicum 2002, 91 (No
2), pp. 263 -
282.\\
20. Cianchi A. Hardy inequalities in Orlicz Spaces. Trans. AMS, {\bf 351},
No 6 (1999), 2456 - 2478.\\
21. Krein S.G., Petunin Yu., and Semenov E.M. Interpolation of linear
operators. AMS, 1982.\\

\newpage

\vspace{14mm}

\begin{center}
Ostrovsky Eugene, Sirota Leonid.\\

Moment Banach Spaces: theory and applications.\\

Abstract.\\
\end{center}
 In this article we introduce and investigate a new class of Banach spaces,
so - called moment spaces, i.e.  which are based on the classical $ L(p) $ spaces,
study their properties: separability, reflexivity, embedding theorems etc.,
and describe some applications to the theory of Fourier series and transform,
theory of martingales, and singular integral operators.\par
 Key words: Banach spaces, moments, Fourier series and transform, martingales,
singular operators.\par
 References: 21 works. \\

\end{document}